\documentclass[draft]{amsart}

\usepackage{mathrsfs}
\usepackage{syntonly}
\usepackage{amsmath}
\usepackage{amsthm}
\usepackage{amsfonts}
\usepackage{amssymb}
\usepackage{latexsym}
\usepackage[colorlinks=true,linkcolor=red,citecolor=blue]{hyperref}
\usepackage[all]{xy}
\usepackage{enumitem}

\newtheorem{intthm}{Theorem}[]
\newtheorem{intcor}{Corollary}[]

\newcommand{\numberseries}{\bfseries}   
\newlength{\thmtopspace}                
\newlength{\thmbotspace}                
\newlength{\thmheadspace}               
\newlength{\thmindent}                  
\newtheoremstyle{bfupright head,slanted body}
                {\thmtopspace}{\thmbotspace}
                {\slshape}{\thmindent}{\bfseries}{.}{\thmheadspace}
                {{\numberseries \thmnumber{#2\;}}\thmnote{#3}}

\newtheoremstyle{bfupright head,upright body}
                {\thmtopspace}{\thmbotspace}
                {\upshape}{\thmindent}{\bfseries}{.}{\thmheadspace}
                {{\numberseries \thmnumber{#2\;}}\thmnote{#3}}

\newtheoremstyle{fixed bf head,slanted body}
                {\thmtopspace}{\thmbotspace}{\slshape}
                {\thmindent}{\bfseries}{.}{\thmheadspace}
                {{\numberseries \thmnumber{#2\;}}\thmname{#1}\thmnote{ (#3)}}

\newtheoremstyle{fixed bf head,upright body}
                {\thmtopspace}{\thmbotspace}{\upshape}
                {\thmindent}{\bfseries}{.}{\thmheadspace}
                {{\numberseries \thmnumber{#2\;}}\thmname{#1}\thmnote{ (#3)}}

\newtheoremstyle{numbered paragraph}
                {\thmtopspace}{\thmbotspace}{\upshape}
                {\thmindent}{\upshape}{}{\thmheadspace}
                {{\numberseries \thmnumber{#2.}}}

\theoremstyle{bfupright head,slanted body}
\newtheorem{res}{}[section]             \newtheorem*{res*}{}

\theoremstyle{bfupright head,upright body}
\newtheorem{bfhpg}[res]{}               \newtheorem*{bfhpg*}{}

\theoremstyle{fixed bf head,slanted body}
\newtheorem{thm}[res]{Theorem}          \newtheorem*{thm*}{Theorem}
\newtheorem{prp}[res]{Proposition}      \newtheorem*{prp*}{Proposition}
\newtheorem{cor}[res]{Corollary}        \newtheorem*{cor*}{Corollary}
\newtheorem{lem}[res]{Lemma}            \newtheorem*{lem*}{Lemma}

\theoremstyle{fixed bf head,upright body}
\newtheorem{dfn}[res]{Definition}       \newtheorem*{dfn*}{Definition}
\newtheorem{rmk}[res]{Remark}           \newtheorem*{rmk*}{Remark}
\newtheorem{exa}[res]{Example}           \newtheorem*{exa*}{Example}
           \newtheorem*{que*}{Question}
           \newtheorem*{fact*}{Fact}
\newtheorem{setup}[res]{Setup}           \newtheorem*{setup*}{Setup}

  \newtheorem*{notation*}{Notation}
  \newtheorem*{complexes*}{Complexes}
  \newtheorem*{morphisms*}{Morphisms}

\theoremstyle{numbered paragraph}
\newtheorem{ipg}[res]{}

\newlength{\thmlistleft}        
\newlength{\thmlistright}       
\newlength{\thmlistpartopsep}   
\newlength{\thmlisttopsep}      
\newlength{\thmlistparsep}      
\newlength{\thmlistitemsep}     

\newcounter{eqc}
\newenvironment{eqc}{\begin{list}{\upshape (\textit{\roman{eqc}})}%
    {\usecounter{eqc}%
      \setlength{\leftmargin}{\thmlistleft}%
      \setlength{\labelwidth}{\thmlistleft}%
      \setlength{\rightmargin}{\thmlistright}%
      \setlength{\partopsep}{\thmlistpartopsep}%
      \setlength{\topsep}{\thmlisttopsep}%
      \setlength{\parsep}{\thmlistparsep}%
      \setlength{\itemsep}{\thmlistitemsep}}}%
  {\end{list}}%

\newcounter{prt}
\newenvironment{prt}{\begin{list}{\upshape (\alph{prt})}%
    {\usecounter{prt}%
      \setlength{\leftmargin}{\thmlistleft}%
      \setlength{\labelwidth}{\thmlistleft}%
      \setlength{\rightmargin}{\thmlistright}%
      \setlength{\partopsep}{\thmlistpartopsep}%
      \setlength{\topsep}{\thmlisttopsep}%
      \setlength{\parsep}{\thmlistparsep}%
      \setlength{\itemsep}{\thmlistitemsep}}}%
  {\end{list}}%

\newcounter{rqm}
  {\end{list}}%

\newenvironment{prf*}[1][Proof]{%
  \begin{proof}[\bf #1]
    \setcounter{equation}{0}
    }
  {\end{proof}
}

\newcommand{\pgref}[1]{\ref{#1}}

\renewcommand{\eqref}[1]{(\pgref{eq:#1})}

\newcommand{\eqclbl}[1]{{\upshape(\textit{#1})}}
\newcommand{\proofofimp}[3][:]{\mbox{\eqclbl{#2}$\!\implies\!$\eqclbl{#3}#1}}

\numberwithin{equation}{res}
\def\urltilda{\kern -.15em\lower .7ex\hbox{\~{}}\kern .04em}

\newcommand{\Rop}{R^\circ}
\newcommand{\oo}{\otimes}
\newcommand{\ra}{\rightarrow}
\newcommand{\Ext}{\mbox{\rm Ext}}
\newcommand{\Hom}{\mbox{\rm Hom}}

\newcommand{\id}{\mbox{\rm id}}
\newcommand{\fd}{\mbox{\rm fd}}
\newcommand{\coker}{\mbox{\rm Coker}}
\newcommand{\ke}{\mbox{\rm Ker}}

\newcommand{\Gfc}{\mathrm{Gfc}\text{-}\mathrm{gldim}}
\newcommand{\Gw}{\mathrm{Gwgldim}}
\newcommand{\Gg}{\mathrm{Ggldim}}
\newcommand{\FFD}{\mathrm{FFD}}
\newcommand{\Gfd}{\mathrm{Gfd}}
\newcommand{\AMod}{A\text{-}\mathsf{Mod}}

\newcommand{\LMod}{\Lambda\text{-}\mathsf{Mod}}

\newcommand{\BMod}{B\text{-}\mathsf{Mod}}
\newcommand{\ModA}{\mathsf{Mod}\text{-}A}
\newcommand{\ModB}{\mathsf{Mod}\text{-}B}
\newcommand{\ModR}{\mathsf{Mod}\text{-}R}

\newcommand{\ABMod}{A\times B\text{-}\mathsf{Mod}}
\newcommand{\RMod}{R\text{-}\mathsf{Mod}}

\newcommand{\eMod}{(B\ltimes_{\theta}M)\text{-}\mathsf{Mod}}

\newcommand{\Mode}{\mathsf{Mod}\text{-}(B\ltimes_{\theta}M)}

\newcommand{\G}{B\ltimes_{\theta}M}
\newcommand{\dg}{\mathsf{dg}}

\begin{document}

\title[Homological invariant properties under cleft extensions]{Homological invariant properties under cleft extensions}%

\author[L. Liang]{Li Liang}

\address{L. Liang: Department of Mathematics, Gansu Center for Fundamental Research in Complex Systems Analysis and Control, Lanzhou Jiaotong University, Lanzhou 730070, China; and Gansu Provincial Research Center for Basic Disciplines of Mathematics and Statistics, Lanzhou 730070, China}

\email{lliangnju@gmail.com}

\urladdr{https://sites.google.com/site/lliangnju}

\author[Y.J. Ma]{Yajun Ma$^{\ast}$}

\address{Y.J. Ma: Department of Mathematics, Gansu Center for Fundamental Research in Complex Systems Analysis and Control, Lanzhou Jiaotong University, Lanzhou 730070, China}

\email{13919042158@163.com}

\thanks{$^{\ast}$ Corresponding author.}

\author[G. Yang]{Gang Yang}

\address{G. Yang: Department of Mathematics, Gansu Center for Fundamental Research in Complex Systems Analysis and Control, Lanzhou Jiaotong University, Lanzhou 730070, China}

\email{yanggang@mail.lzjtu.cn}





\begin{abstract}
We study the behavior of the Gorenstein weak global dimension under a cleft extension of rings; we prove that under some mild conditons the finiteness of the Gorenstein weak global dimension is invariant. Moreover, we compare the relative singularity categories with respect to flat-cotorsion modules under a cleft extension of rings. Some applications to $\theta$-extensions and Morita context rings are given.
\end{abstract}

\maketitle

\thispagestyle{empty}

\section{Introduction}
\noindent
As one of the most basic homological invariants attached to an associative ring, the weak global dimension defined in terms of the flat dimension is left-right symmetric. In Gorenstein homological algebra, there are corresponding definition and conclusion: The Gorenstein weak global dimension defined in terms of the Gorenstein flat dimension is also left-right symmetric (see \cite[Corollary 2.5]{CET}). As a consequence, one gets a connection between the rings of finite Gorenstein weak global dimension and Ding-Chen rings; it was proved in \cite[Proposition 3.4]{CELTW} that a coherent ring $R$ is Ding-Chen if and only if $R$ has finite Gorenstein weak global dimension.

Our first aim in the paper is to study the behavior of the Gorenstein weak global dimension under a cleft extension of rings, which is a ring homomorphism $f: A\to B$ together with a ring homomorphism $g: B\to A$ such that $fg=\mathrm{Id}_{B}$. We mention that a cleft extension $f: A\to B$ of rings gives rise a diagram
\begin{equation}\label{cleft}
  \xymatrix@C=3pc{
   \BMod\ar[r]^{\sf i} & \AMod
   \ar[r]^{\sf e} & \BMod
    \ar@/_1.8pc/[l]_{{\sf l}}
  }
\end{equation}
such that $\mathsf{ei}\simeq\mathrm{Id}_{B}$, where $\sf i$ and $\sf e$ are the restriction functors along $f$ and $g$, respectively, and ${\sf l}=A\otimes_B-$. The above diagram is called a cleft extension of module categories, and is denoted by $(\BMod,\AMod,\mathsf{i},\mathsf{e},\mathsf{l})$.

We mention that every cleft extension of module categories occurs as a cleft extension associated to a $\theta$-extension (see Remark \ref{remark:fg}). So we can summarize the first main result on the Gorenstein weak global dimension in the context of $\theta$-extensions to make it clearer to readers, which is listed in Theorem \ref{IT:corollary}.

\begin{intthm}\label{thmA}
Let $M$ be a $B$-bimodule with $M^{\oo s}=0$ for some positive integer $s$, and let $\theta:M\oo_{B}M\ra M$ be an associative $B$-bimodule homomorphism. Consider the $\theta$-extension $\G$. If $M$ is projective as a left $B$-module and a right $B$-module, then there are inequalities
    $$\Gw(B)\leq\Gw(\G)\leq \Gw(B)+s-1.$$
    In particular, $\Gw(B)$ is finite if and only if $\Gw(\G)$ is finite.
\end{intthm}

\begin{rmk*}
There is another important homological invariant attached to an associative ring $R$: the Gorenstein global dimension $\Gg(R)$, which is defined in terms of the Gorenstein projective dimension or the Gorenstein injective dimension. It is known that there is an inequality $\Gw(R)\leq\Gg(R)$ by \cite[Theorem 5.13]{CELTW}; the equality is not true in general. By the proof of \cite[Propositions 4.3 and 4.4]{Kostas} and using \cite[Theorem 4.1]{Emmanouil}, one gets that the inequalities
$$\Gg(B)\leq\Gg(\G)\leq \Gg(B)+s-1$$
hold under the same assumptions as in Theorem \ref{thmA}.
\end{rmk*}

We mention that the Morita context rings with zero bimodule homomorphisms can be viewed as a trivial extension, and consequently as a $\theta$-extension; see \cite[Proposition 2.5]{GP}.
Hence as an application of Theorem \ref{thmA}, we have the following result.

\begin{intcor}
Let $\Lambda =\left(\begin{smallmatrix}  A & {_{A}}N_{B}\\  {_{B}}M_{A} & B \\\end{smallmatrix}\right)$
be a Morita context ring with $N\oo_{B}M=0=M\oo_{A}N$.
If $M$ is projective as a right $A$-module and a left $B$-module, and $N$ is projective as a left $A$-module and a right $B$-module, then there are inequalities
$$\max\{\Gw (A),\Gw (B)\}\leq\Gw (\Lambda)$$
and
$$\Gw (\Lambda)\leq \max\{\Gw (A),\Gw (B)\}+1.$$
In particular, $\Gw(\Lambda)$ is finite if and only if the invariants $\Gw(A)$ and $\Gw(B)$ are finite.
\end{intcor}

A cleft extension $(\BMod,\AMod,\mathsf{i},\mathsf{e},\mathsf{l})$ of module categories gives rise to an endofunctor $\mathsf{F}: \BMod\ra \BMod$, which carries significant implications. Our second aim in the paper is to explore the relative singularity categories with respect to flat-cotorsion modules in the context of a cleft extension of module categories with the induced endofunctor $\sf F$ both nilpotent and f-perfect (see Theorem \ref{thm:relative singularity categories}). Here, the key tool we use is the flat model structure on the category of chain complexes constructed by Gillespie, which offers a different perspective compared to the Kostas's work \cite{Kostas}.

\begin{intthm}
Let $A$ and $B$ be right coherent rings, and let $(\BMod,\AMod,\mathsf{i},\mathsf{e},\mathsf{l})$ be a cleft extension of module categories (see (\ref{cleft})) such that the induced endofunctor $\mathsf{F}: \BMod\to \BMod$ is f-perfect and nilpotent. Assume that the following conditions are satisfied:
\begin{prt}
\item $\sf l$  and $\sf r$ are exact, where $\sf r$ is the right adjoint functor of $\sf e$.
\item $\sf l$ and $\sf q$ preserve products, where $\sf q$ is the left adjoint functor of $\sf i$.
\item $\sf i$ preserves cotorsion modules.
\end{prt}
Then there is a diagram
\begin{equation*}
  \xymatrix@C=2pc{
   \mathbf{D}_{\mathsf{FlatCot}}(B)\ar[r]^{\overline{\sf Ri}} & \mathbf{D}_{\mathsf{FlatCot}}(A) \ar@/_1.8pc/[l]_{\overline{\mathsf{Lq}}}
   \ar[r]^{\overline{\sf Re}} & \mathbf{D}_{\mathsf{FlatCot}}(B) \ar@/_1.8pc/[l]_{\overline{\mathsf{Ll}}}}
\end{equation*}
of relative singularity categories and triangle functors such that $(\overline{\sf Lq},\overline{\sf Ri})$ and $(\overline{\sf Ll},\overline{\sf Re})$ are adjoint pairs, $\overline{\sf Re}\circ\overline{\sf Ri}\simeq \mathrm{Id}_{\mathbf{D}_{\mathsf{FlatCot}}(B)}$ and $\overline{\sf Lq}\circ\overline{\sf Ll}\simeq \mathrm{Id}_{\mathbf{D}_{\mathsf{FlatCot}}(B)}$.
\end{intthm}

An application to Morita context rings can be found in Corollary \ref{cor6.6}.

\begin{intcor}
Let $\Lambda =\left(\begin{smallmatrix}  A & {_{A}}N_{B}\\  {_{B}}M_{A} & B \\\end{smallmatrix}\right)$
be a Morita context ring with $A$ and $B$ right coherent and $N\oo_{B}M=0=M\oo_{A}N$.
Assume that $M$ is projective as a right $A$-module and a left $B$-module, and $N$ is projective as a left $A$-module and a right $B$-module, and assume that $M$ and $N$ are finitely generated as a right $A$-module and a right $B$-module respectively. Then there is a diagram
\begin{equation*}
  \xymatrix@C=2pc{
   \mathbf{D}_{\mathsf{FlatCot}}(A)\times\mathbf{D}_{\mathsf{FlatCot}}(B) \ar[r]^-{\overline{\sf Ri}} & \mathbf{D}_{\mathsf{FlatCot}}(\Lambda) \ar@/_1.8pc/[l]_{\overline{\mathsf{Lq}}}
   \ar[r]^-{\overline{\sf Re}} & \mathbf{D}_{\mathsf{FlatCot}}(A)\times\mathbf{D}_{\mathsf{FlatCot}}(B) \ar@/_1.8pc/[l]_-{\overline{\mathsf{Ll}}}}
\end{equation*}
of relative singularity categories and triangle functors such that $(\overline{\sf Lq},\overline{\sf Ri})$ and $(\sf \overline{Ll},\overline{Re})$ are adjoint pairs, $\overline{\sf Re}\circ\overline{\sf Ri}\simeq \mathrm{Id}_{\mathbf{D}_{\mathsf{FlatCot}}(A)\times\mathbf{D}_{\mathsf{FlatCot}}(B)}$ and $ \overline{\sf Lq}\circ\overline{\sf Ll}\simeq \mathrm{Id}_{\mathbf{D}_{\mathsf{FlatCot}}(A)\times\mathbf{D}_{\mathsf{FlatCot}}(B)}$.
\end{intcor}

The paper is organized as follows. In Section 2, we recall some notations, the definitions of cleft extensions and cleft coextensions of abelian categories, together with their basic homological properties required in proofs. In Section 3, we study the Gorenstein weak global dimension under a cleft extension of module categories.
In Section 4, our main results are applied to some known ring extensions.
In Section 5, we compare the relative singularity categories of rings linked by cleft extensions of module categories and apply the results to the Morita context rings.

\section{Preliminaries on cleft extensions}
\noindent

In this section, we fix some notation, recall relevant notions and collect some necessary facts on cleft extensions.

\begin{setup}
All rings considered in the paper are nonzero associative rings with identity and all modules are unitary. For a ring $R$, we usually work with left $R$-modules and write $\RMod$ for the category of left $R$-modules. Similarly, we write $\ModR$ for the category of right $R$-modules.
We use $_{R}M$ (resp., $M_{R}$) to denote a left (resp., right) $R$-module $M,$ and denote the projective dimension, injective dimension and flat dimension of $_{R}M$ (resp., $M_{R}$) by $\mathrm{pd}_{R}M, \mathrm{id}_{R}M$ and $\mathrm{fd}_{R}M$ (resp., $\mathrm{pd}M_{R}, \mathrm{id}M_{R}$ and $\mathrm{fd}M_{R}$).
\end{setup}

\begin{bfhpg}[Cleft extension]\label{se.F}
Recall from Beligiannis \cite{Be1} that a {\it cleft extension} of an abelian category $\mathsf{B}$ is an abelian category $\mathsf{A}$ together with functors
\begin{equation*}
  \xymatrix@C=3pc{
   \mathsf{B} \ar[r]^-{\sf i} & \mathsf{A}
   \ar[r]^{\sf e} & \mathsf{B} \ar@/_1.8pc/[l]_{\sf l}}
\end{equation*}
satisfying the following conditions:
\begin{prt}
\item The functor $\mathsf{e}$ is faithful and exact.
\item The pair $(\mathsf{l},\mathsf{e})$ is an adjoint pair.
\item There is a natural isomorphism $\varphi: \mathsf{ei}\ra \mathrm{Id}_{\mathsf{B}}$.
\end{prt}

We always denote the above cleft extension by $(\mathsf{B},\mathsf{A},\mathsf{i},\mathsf{e},\mathsf{l})$;
the data can give rise to more structure information.
For instance, it follows from \cite[Lemma 2.2]{GP1} that $\mathsf{i}$ is fully faithful and exact. Moreover, there is a functor $\mathsf{q}:\mathsf{A}\to \mathsf{B}$ such that $(\mathsf{q},\mathsf{i})$ forms an adjoint pair.

We will give some endofunctors on the categories $\mathsf{A}$ and $\mathsf{B}$, which are important for our results. Denote by $\nu:\rm Id_{\mathsf{B}}\to \mathsf{el}$ the unit and by $\mu: \mathsf{le}\to \mathrm{Id}_{\mathsf{A}}$ the counit of the adjoint pair $(\mathsf{l},\mathsf{e})$.
The unit and counit satisfy
$$\mathrm{Id}_{\mathsf{l}(B)}=\mu_{\mathsf{l}(B)}\mathsf{l}(\nu_{B})\ \ \mathrm{and}\ \
\mathrm{Id}_{\mathsf{e}(A)}=\mathsf{e}(\mu_{A})\nu_{\mathsf{e}(A)}$$
for all $A\in\mathsf{A}$ and $B\in\mathsf{B}$.
It follows that $\mathsf{e}(\mu_{A})$ is a split epimorphism. Since $\sf e$ is faithful exact, it implies that $\mu_{A}: \mathsf{le}(A)\ra A$ is an epimorphism.
Hence we have the following exact sequence
$$0\ra \ke \mu_{A}\to \mathsf{le}(A)\xrightarrow{\mu_{A}} A\ra 0$$
for every $A\in \mathsf{A}.$
The assignment $A\mapsto\ke \mu_{A}$ defines an endofunctor $\mathsf{G}:\mathsf{A}\to \mathsf{A}$.
Fix an object $B$ in $\mathsf{B}$.
Then we have a short exact sequence in $\mathsf{A}$:
$$0\ra \mathsf{G}(\mathsf{i}(B))\ra \mathsf{le}(\mathsf{i}(B))\xrightarrow{\mu_{\mathsf{i}(B)}} \mathsf{i}(B)\ra 0.$$
We denote by $\mathsf{F}(B)=\mathsf{e}(\mathsf{G}(\mathsf{i}(B))).$
Viewing the isomorphism $\mathsf{ei}(B)\cong B$ as an identification and applying the exact functor $\mathsf{e}:\mathsf{A}\to \mathsf{B}$ we get the exact sequence:
\begin{equation}
\label{se.F}
0\ra \mathsf{F}(B)\ra \mathsf{el}(B)\xrightarrow{\mathsf{e}(\mu_{\mathsf{i}(B)})}B\ra 0.
\end{equation}
The assignment $B\mapsto \mathsf{F}(B)$ defines an endofunctor $\mathsf{F}:\mathsf{B}\to \mathsf{B}.$
Since $\mathsf{e}(\mu_{\mathsf{i}(B)})$ is a split epimorphism, it follows that there is natural isomorphism $\mathsf{el}\simeq \mathsf{F}\oplus \mathrm{Id}_{\sf B}.$
\end{bfhpg}

We collect some results from \cite{GP1,Kostas} in the following, which will be used later.

\begin{lem}\label{lem of cleft extension}
Let $(\mathsf{B},\mathsf{A},\mathsf{i},\mathsf{e},\mathsf{l})$ be a cleft extension of abelian categories with $\sf F: \mathsf{B}\to \mathsf{B}$ and $\sf G: \mathsf{A}\to \mathsf{A}$ the induced endofunctors. Then the following statements hold:
\begin{prt}
\item For each $n\geq 1$, one has $\mathsf{F}^{n}=0$ if and only if $\mathsf{G}^{n}=0$.
\item For each $n\geq 1$ and every $A\in \mathsf{A}$, there is an exact sequence
    $$0\to \mathsf{G}^{n}(A)\to \mathsf{l}\mathsf{F}^{n-1}\mathsf{e}(A)\to \mathsf{G}^{n-1}(A)\to 0.$$
\end{prt}
\end{lem}

\begin{bfhpg}[Cleft coextension]
Recall from \cite[Definition 2.1]{Be1} that a {\it cleft coextension} of an abelian category $\mathsf{B}$ is an abelian category $\mathsf{A}$ together with functors
\begin{equation*}
  \xymatrix@C=3pc{
    \mathsf{B} \ar[r]^-{\sf i} & \mathsf{A} \ar[r]^-{\sf e} & \mathsf{B}
    \ar@/^1.8pc/[l]^-{\sf r}}
\end{equation*}
satisfying the following conditions:
\begin{prt}
\item The functor $\sf e$ is faithful exact.
\item The pair $(\sf e,\sf r)$ is an adjoint pair.
\item There is a natural isomorphism $\varphi:{\sf ei}\ra \mathrm{Id}_{\mathsf{B}.}$
\end{prt}

We always denote a cleft coextension by $(\mathsf{B},\mathsf{A},\mathsf{i},\mathsf{e},\mathsf{r})$. We will give some endofunctors on the categories $\mathsf{A}$ and $\mathsf{B}$. Denote by $\nu':\rm Id_{\mathsf{A}}\ra \sf re$ the unit and by $\mu':{\sf er}\ra \mathrm{Id}_{\mathsf{B}}$ the counit of the adjoint pair $(\sf e,r)$.
The unit and counit satisfy
$$ \mathrm{Id}_{\mathsf{e}(A)}=\mu'_{\mathsf{e}(A)}\mathsf{e}(\nu'_{A})\ \  \mathrm{and}\ \  \mathrm{Id}_{\mathsf{r}(B)}=\mathsf{r}(\mu'_{B})\nu'_{\mathsf{r}(B)}$$
for all $A\in\mathsf{A}$ and $B\in\mathsf{B}$.
It follows that $\mathsf{e}(\nu'_{A})$ is a (split) monomorphism. Since $\sf e$ is faithful exact, we have $\nu'_{A}: A\ra \mathsf{re}(A)$ is a monomorphism.
Hence we have the following exact sequence
$$0\ra A\xrightarrow{\nu'_{A}} \mathsf{re}(A)\ra \coker\nu'_{A} \ra 0$$
for every $A\in \mathsf{A}.$
The assignment $A\mapsto\coker \nu'_{A}$ defines an endofunctor $\mathsf{G}':\mathsf{A}\ra \mathsf{A}$.
Consider now an object $B$ in $\mathsf{B}$.
Then we have a short exact sequence in $\mathsf{A}$:
$$0\ra \mathsf{i}(B)\xrightarrow{\nu'_{\mathsf{i}(B)}}\mathsf{re}(\mathsf{i}(B)) \ra \mathsf{G}'(\mathsf{i}(B)) \ra 0.$$
We denote by $\mathsf{F}'(B)=\mathsf{e}(\mathsf{G}'(\mathsf{i}(B))).$
Viewing the isomorphism $\mathsf{ei}(B)\cong B$ as an identification and applying the exact functor $\mathsf{e}:\mathsf{A}\ra \mathsf{B}$ we get the exact sequence:
$$0\ra B\xrightarrow{\mathsf{e}(\nu'_{\mathsf{i}(B)})} \mathsf{er}(B)\ra \mathsf{F}'(B)\ra 0.$$
The assignment $B\mapsto \mathsf{F}'(B)$ defines an endofunctor $\mathsf{F}':\mathsf{B}\ra \mathsf{B}.$
Since ${\sf e}(\nu'_{\mathsf{i}(B)})$ is a split monomorphism, it follows that there is a natural isomorphism $\mathsf{er}\simeq \mathsf{F}'\oplus \mathrm{Id}_{\sf B}.$
\end{bfhpg}

\begin{rmk}\label{adjoint}
We mention that a cleft extension $(\mathsf{B},\mathsf{A},\mathsf{i},\mathsf{e},\mathsf{l})$ is the upper part of a cleft coextension $(\mathsf{B},\mathsf{A},\mathsf{i},\mathsf{e},\mathsf{r})$ if and only if the induced endofunctor $\mathsf{F}$ appearing in (\ref{se.F}) has a right adjoint $\mathsf{F}'$. Here $\mathsf{F}'$ is actually the induced endofunctor by the cleft coextension $(\mathsf{B},\mathsf{A},\mathsf{i},\mathsf{e},\mathsf{r})$; see \cite[Remark 2.7]{Be1}.
\end{rmk}

The next result is a dual version of Lemma \ref{lem of cleft extension}.

\begin{lem}\label{lem of cleft coextension}
Let $(\mathsf{B},\mathsf{A},\mathsf{i},\mathsf{e},\mathsf{r})$ be a cleft coextension of abelian categories with $\sf F': \mathsf{B}\to \mathsf{B}$ and $\sf G': \mathsf{A}\to \mathsf{A}$ the induced endofunctors. Then the following statements hold:
\begin{prt}
\item Let each $n\geq 1$, one has $\mathsf{F}'^{n}=0$ if and only if $\mathsf{G}'^{n}=0$.
\item For each $n\geq 1$ and every $A\in \mathsf{A}$, there is an exact sequence
    $$0\ra \mathsf{G}'^{n-1}(A)\ra \mathsf{r}\mathsf{F}'^{n-1}\mathsf{e}(A)\ra \mathsf{G}'^{n}(A)\ra 0.$$
\end{prt}
\end{lem}

We now turn our attention to cleft extensions of module categories. The next result is from \cite[Proposition 4.19]{Kostas}, which will be used frequently. {\sf In the rest of the paper}, we let $A$ and $B$ two associative rings. Recall that an endofunctor $\sf F$ on $\BMod$ is called {\it nilpotent} if there is an integer $n\geq 1$ such that ${\sf F}^n=0$.

\begin{prp}\label{extension-coextension}
A cleft extension $(\BMod,\AMod,\mathsf{i},\mathsf{e},\mathsf{l})$ of module categories with $\sf F$ the induced endofunctor on $\BMod$ is the upper part of a cleft coextension $(\BMod,\AMod,\mathsf{i},\mathsf{e},\mathsf{r})$ of module categories with $\sf F'$ the induced endofunctor on $\BMod$, that is, there exists a diagram
\begin{equation}\label{2.7.1}
  \xymatrix@C=2pc{
    \BMod \ar[r]^-{\sf i} & \AMod \ar[r]^-{\sf e} & \BMod.
    \ar@/_1.6pc/[l]_-{\sf l}
    \ar@/^1.6pc/[l]^-{\sf r}\ar@(lu,ru)[]^{\sf F}\ar@(ld,rd)[]_{\mathsf{F}'}}
\end{equation}
Moreover, if $\mathsf{F}$ is  nilpotent, then so is $\mathsf{F}'$.
\end{prp}


\begin{rmk}
We mention that $(\mathsf{F},\mathsf{F}')$ forms an adjoint pair; see Remark \ref{adjoint}. Thus $\mathsf{F}$ is a right exact functor and $\mathsf{F}'$ is a left exact functor.
\end{rmk}

In the following we give some examples of cleft extensions of module categories which will be used later in the paper.

\begin{exa}\label{ex:extension}
Let $M$ be a $B$-bimodule, and let $\theta:M\oo_{B}M\ra M$ be an associative $B$-bimodule homomorphism, that is, $\theta$ is a bimodule homomorphism such that $\theta(\theta\otimes \mathrm{Id}_{M})=\theta(\mathrm{Id}_{M}\otimes \theta).$
Then $\theta$-\emph{extension} of $B$ by $M$, denoted by $B\ltimes_{\theta}M$, is defined to be the ring with underlying group $B\oplus M$ and the multiplication given as follows:
$$(b,m)\cdot(b',m')=(bb',b m'+mb'+\theta(m\oo m')),$$
for $b,b'\in B$ and $m,m'\in M$.

Consider the ring homomorpisms $f:B\ltimes_{\theta}M\ra B$ and $g:B\ra B\ltimes_{\theta}M$ given by $(b,m)\mapsto b$ and $b\mapsto(b,0)$, respectively.
Denote by $\mathsf{Z}_{1}:\BMod\ra \eMod$ and $\mathsf{Z}_{2}:\ModB\ra \Mode$  the
restriction functors induced by $f$, and denoted by $\mathsf{U}_{1}:\eMod\ra \BMod$ and
$\mathsf{U}_{2}:\Mode\ra \ModB$ the restriction functors induced by $g$.
Then there are (co)cleft extensions of module categories
\begin{equation}\label{2.8.1}
  \xymatrix@C=2pc{
    \BMod \ar[r]^-{\mathsf{Z}_{1}} & \eMod \ar[r]^-{\mathsf{U}_{1}} \ar@/_1.8pc/[l]_-{\mathsf{C}_{1}=B\oo_{{B\ltimes_{\theta}M}}-}
    \ar@/^1.8pc/[l]^-{\mathsf{K}_{1}=\mathrm{Hom}_{B\ltimes_{\theta}M}(B,-)}
     & \BMod
    \ar@/_1.8pc/[l]_-{\mathsf{T}_{1}=(B\ltimes_{\theta}M)\oo_{B}-}
    \ar@/^1.8pc/[l]^-{\mathsf{H}_{1}=\mathrm{Hom}_{B}(B\ltimes_{\theta}M,-)}}
\end{equation}
and
\begin{equation}\label{2.8.2}
  \xymatrix@C=2pc{
    \ModB \ar[r]^-{\mathsf{Z}_{2}} & \Mode \ar[r]^-{\mathsf{U}_{2}}
    \ar@/_1.8pc/[l]_-{\mathsf{C}_{2}=-\oo_{B\ltimes_{\theta}M}B}
    \ar@/^1.8pc/[l]^-{\mathsf{K}_{2}=\mathrm{Hom}_{B\ltimes_{\theta}M}(B,-)}
    & \ModB.
    \ar@/_1.8pc/[l]_-{\mathsf{T}_{2}=-\oo_{B}(B\ltimes_{\theta}M)}
    \ar@/^1.8pc/[l]^-{\mathsf{H}_{2}=\mathrm{Hom}_{B}(B\ltimes_{\theta}M,-)}}
\end{equation}
We mention that the induced endofunctor $\mathsf{F}_{1}$ (resp., $\mathsf{F}_{2}$) on $\BMod$ (resp., $\ModB$) is naturally isomorphic to $M\oo_{B}-$ (resp., $-\oo_{B}M$).
\end{exa}

\begin{rmk} \label{remark:fg}
We mention that every cleft extension of module categories occurs as a cleft extension associated to a $\theta$-extension. More details, given a a cleft extension of module categories $(\BMod,\AMod,\mathsf{i}_{1},\mathsf{e}_{1},\mathsf{l}_{1})$, it follows from \cite[Proposition 6.9]{KP} that there are a $B$-bimodule $M$ and an associative $B$-bimodule homomorphism $\theta:M\oo_{B}M\ra M$ such that the following diagram commutes:
\begin{equation*}
  \xymatrix@C=3pc{
    \BMod \ar@{=}[d] \ar[r]^-{\mathsf{i}_{1}} & \AMod \ar[d]_-{\simeq} \ar@/_1.8pc/[l]_-{\mathsf{q}_{1}} \ar[r]^-{\mathsf{e}_{1}} & \BMod \ar@{=}[d]
    \ar@/_1.8pc/[l]_-{\mathsf{l}_{1}}\\
    \BMod \ar[r]^-{\mathsf{Z}_{1}} & \eMod \ar@/_1.8pc/[l]_-{\mathsf{C}_{1}} \ar[r]^-{\mathsf{U}_{1}} & \BMod
    \ar@/_1.8pc/[l]_-{\mathsf{T}_{1}}}
\end{equation*}
where $\AMod\xrightarrow{\simeq} \eMod$ is an equivalence, $\mathsf{q}_{1}$ is the left adjoint of $\mathsf{i}_{1}$, and the bottom diagram is the cleft extension associated to $B\ltimes_{\theta}M$; see (\ref{2.8.1}).

Consider the $\theta$-extension $B\ltimes_{\theta}M$. Then there is a cleft extension
\begin{equation*}
  \xymatrix@C=2pc{
    \ModB \ar[r]^-{\mathsf{Z}_{2}} & \Mode \ar[r]^-{\mathsf{U}_{2}}
    \ar@/_1.8pc/[l]_-{\mathsf{C}_{2}=-\oo_{B\ltimes_{\theta}M}B}
    & \ModB
    \ar@/_1.8pc/[l]_-{\mathsf{T}_{2}=-\oo_{B}(B\ltimes_{\theta}M)}}
\end{equation*}
of module categories; see (\ref{2.8.2}). It follows from \cite[Proposition 18.32]{Lam} that there is an equivalence $\ModA\xrightarrow{\simeq} \Mode$, and so one gets a cleft extension $(\ModB,\ModA,\mathsf{i}_{2},\mathsf{e}_{2},\mathsf{l}_{2})$ such that the following diagram commutes:
\begin{equation*}
  \xymatrix@C=3pc{
    \ModB \ar@{=}[d] \ar[r]^-{\mathsf{i}_{2}} & \ModA \ar[d]_-{\simeq} \ar@/_1.8pc/[l]_-{\mathsf{q}_{2}} \ar[r]^-{\mathsf{e}_{2}} & \ModB \ar@{=}[d]
    \ar@/_1.8pc/[l]_-{\mathsf{l}_{2}}\\
    \ModB \ar[r]^-{\mathsf{Z}_{2}} & \Mode \ar@/_1.8pc/[l]_-{\mathsf{C}_{2}} \ar[r]^-{\mathsf{U}_{2}} & \ModB,
    \ar@/_1.8pc/[l]_-{\mathsf{T}_{2}}}
\end{equation*}
where $\mathsf{q}_{2}$ is the left adjoint of $\mathsf{i}_{2}$.
We always denote by $\mathsf{F}_{1}$ and $\mathsf{F}_{2}$ the induced endofunctors on $\BMod$ and $\ModB$, respectively. Then $\mathsf{F}_{1}\simeq M\oo_{B}-$ and $\mathsf{F}_{2}\simeq -\oo_{B}M$.
\end{rmk}

\begin{lem}\label{lem:left-right}
Keep the notation in Remark \ref{remark:fg}. Let $\mathsf{r}_{1}$ be the right adjoint of $\mathsf{e}_{1}$ and $\mathsf{r}_{2}$ the right adjoint of $\mathsf{e}_{2}.$
Then the following statements hold:
\begin{prt}
\item If $\mathsf{r}_{1}$ is exact, then so is $\mathsf{l}_{2}.$
\item If $\mathsf{r}_{2}$ is exact, then so is $\mathsf{l}_{1}.$
\item Let $B$ be a commutative ring. If $\mathsf{r}_{1}$ is exact, then so are $\mathsf{r}_{2},\mathsf{l}_{2}$ and $\mathsf{l}_{1}.$
\end{prt}
\end{lem}
\begin{prf*}
It is enough to consider the case of the cleft extension associated to a $\theta$-extension $B\ltimes_{\theta}M$; see Remark \ref{remark:fg}. We mention that the right adjoint functor  $\mathsf{H}_{1}$ of $\mathsf{U}_{1}$ is exact if and only if $B\ltimes_{\theta}M$ is a projective left $B$-module; see Example \ref{ex:extension}. Then $\mathsf{T}_{2}$ is exact, and so the statement (a) follows. Using a similar method one can prove the statement (b), and get that if $\mathsf{r}_{1}$ is exact then so is $\mathsf{r}_{2}$ under the condition that $B$ is commutative. Then the statement (c) follows from (a) and (b).
\end{prf*}

\begin{exa}\label{ex:Morita ring}
Let $_{B}M_{A}$ and $_{A}N_{B}$ be bimodules.
Consider the Morita context ring $\Lambda =\left(\begin{smallmatrix}  A & {_{A}}N_{B}\\  {_{B}}M_{A} & B \\\end{smallmatrix}\right)$
with zero bimodule homomorphisms.
We mention that a module category over the Morita context ring $\Lambda$ is equivalent to a category whose objects are tuples $(X,Y,f,g)$, where $X\in \AMod$, $Y\in \BMod$, $f\in \text{Hom}_{B}(M\oo_{A}X,Y)$ and $g\in \text{Hom}_{A}(N\oo_{B}Y,X)$; see \cite{eg82}.
Thus there is a cleft (co)extension of module categories
\begin{equation*}
  \xymatrix@C=2pc{
    \ABMod \ar[r]^-{\sf i} & \LMod \ar@/_1.6pc/[l]_-{\sf q} \ar@/^1.6pc/[l]^-{\sf p}
    \ar[r]^-{\sf e} & \ABMod.
    \ar@/_1.6pc/[l]_-{\sf l}
    \ar@/^1.6pc/[l]^-{\sf r}\ar@(lu,ru)[]^{\sf F}\ar@(ld,rd)[]_{\mathsf{F}'}}
\end{equation*}
The functors are given as follows:
\begin{enumerate}
\item $\mathsf{i}(X,Y)=(X,Y,0,0).$
\item $\mathsf{e}(X,Y,f,g)=(X,Y).$
\item $\mathsf{l}(X,Y)=(X,M\oo_{A}X,1,0)\oplus (N\oo_{B}Y,Y,0,1)$.
\item $\mathsf{q}(X,Y,f,g)=(X,\coker f)\oplus (\coker g, Y).$
\item $\mathsf{r}(X,Y)=(X,\Hom_{A}(N,X),0,\varepsilon'_{X})\oplus (\Hom_{B}(M,Y),Y,\varepsilon_{Y},0)$. Here
 $\varepsilon':N\oo_{B}\Hom_{A}(N,-)\ra \mathrm{Id}_{\AMod}$ and $\varepsilon_{}:M\oo_{A}\Hom_{B}(M,-)\ra \mathrm{Id}_{\BMod}$
  are counits of adjoint pairs $(-\oo_{A}N,\Hom_{B}(N,-))$ and $(-\oo_{B}M,\Hom_{A}(M,-))$, respectively.
\item $\mathsf{p}(X,Y,f,g)=(\ke\widetilde{f},\ke\widetilde{g})$. Here $\widetilde{f}:  X\to\Hom_{B}(M,Y)$ is defined as
    $\widetilde{f}(x)(m)=f(m\otimes x)$ for all $x\in X, m\in M$, and
    $\widetilde{g}: Y\to\Hom_{A}(N,X)$ is defined as $\widetilde{g}(y)(n)=g(n\otimes y)$ for all $y\in Y, n\in N$.
\item $\mathsf{F}(X,Y)=(N\oo_{B}Y,M\oo_{A}X)$.
\item $\mathsf{F}'(X,Y)=(\Hom_{B}(N,Y),\Hom_{A}(M,X))$.
\end{enumerate}
\end{exa}

\begin{rmk}
Note that a Morita context ring with zero bimodule homomorphisms can be viewed as a trivial extension; see \cite[Proposition 2.5]{GP}. Thus Example \ref{ex:Morita ring} is a special case of Example \ref{ex:extension}.
\end{rmk}

The next definition of left perfect functors can be found in \cite[Definition 6.4]{KP}.

\begin{dfn}\label{def:left perfect}
A right exact functor $\mathsf{F}:\BMod \ra \BMod$ is called \emph{left perfect} if the following conditions hold:
\begin{prt}
\item $\mathbb{L}_{i}\mathsf{F}^{j}(\mathsf{F}P)=0$ for every projective left $B$-module $P$ and all $i,j\geq 1.$
\item There is an integer $n\geq1$ such that for every $p,q\geq1$ satisfying $p+q\geq n+1,$
one has $\mathbb{L}_{p}\mathsf{F}^{q}=0.$
\end{prt}

A left perfect functor $\mathsf{F}:\BMod \ra \BMod$ is called \emph{f-perfect} if $\mathrm{fd}_{B}\mathsf{F}(P)<\infty$ for every projective left $B$-module $P$; this notion extends the one of perfect functors in the sense of \cite[Definition 3.1]{Kostas}.
\end{dfn}

Recall that a $B$-bimodule $M$ is called \emph{nilpotent} if $M^{\otimes n}=0$ for some integer $n\geq 1$.

\begin{lem}\label{lem:flat dimension}
Let $M$ be a $B$-bimodule with $\mathrm{Tor}_{i}^{B}(M,M^{\oo j})=0$ for all $i,j\geq1$.
If $M_{B}$ has finite flat dimension, then $\mathrm{fd}(M^{\otimes j})_{B}\leq j(\mathrm{fd}M_{B}).$
\end{lem}
\begin{prf*}
Assume that $\mathrm{fd}M_{B}=k,$ and take a flat resolution of $M_{B}$
$$0\to F_{k}\to \cdots \to F_{1}\to F_{0}\to M\to 0.$$
Since $\mathrm{Tor}_{i}^{B}(M,M)=0$ for all $i\geq1$, it follows that the following sequence
$$0\to F_{k}\otimes_{B}M\to \cdots \to F_{1}\otimes_{B}M\to F_{0}\otimes_{B}M\to M\otimes_{B}M\to 0$$
is exact.
We observe that $\mathrm{fd}(F_{i}\otimes_{B}M)_{B}\leq \mathrm{fd}M_{B}$ for all $0\leq i\leq k$.
Then $\mathrm{fd}(M\otimes_{B}M)_{B}\leq 2(\mathrm{fd}M_{B}).$
Inductively, we have that
 $\mathrm{fd}(M^{\otimes j})_{B}\leq j(\mathrm{fd}M_{B}).$
\end{prf*}

\begin{lem}\label{left perfect module}
Let $M$ be a $B$-bimodule with $M^{\otimes s}=0$ for some positive integer $s$. If $M_B$ has finite flat dimension and $\mathrm{Tor}_{i}^{B}(M,M^{\oo j})=0$ for all $i,j\geq1$, then the functor $F=M\otimes_B-$ is left perfect with $\mathbb{L}_{p}\mathsf{F}^{q}=\mathrm{Tor}_{p}^{B}(M^{\otimes q},-)=0$ for every $p,q\geq 1$ satisfying $p+q\geq \sup\{\mathrm{fd}(M^{\otimes j})_{B}\ |\ j\geq1\}+s$.
\end{lem}
\begin{prf*}
By Lemma \ref{lem:flat dimension} we have $\mathrm{fd}(M^{\otimes j})_{B}\leq j(\mathrm{fd}M_{B}).$
Since $M^{\otimes s}=0$ (we may assume that $s\geq2$), it follows that
$\sup\{\mathrm{fd}(M^{\otimes j})_{B}\ |\ j\geq1\}=n_{M}<\infty.$ Hence we have $\mathrm{Tor}_{t}^{B}(M^{\otimes j},-)=0$ for all $j\geq 1$ and $t>n_{M}$, which yields that $\mathrm{Tor}_{p}^{B}(M^{\otimes q},-)=0$ for all $p,q\geq 1$ satisfying $p+q\geq n_{M}+s$ as $M^{\otimes s}=0$. On the other hand, by \cite[Lemma 5.2]{KP} one gets that $\mathrm{Tor}_{i}^{B}(M^{\oo j},M)=0$ for all $i,j\geq1$. Thus the functor $M\otimes_B-$ is left perfect.
\end{prf*}

\begin{prp}\label{perfect module}
Let $M$ be a $B$-bimodule. Then the following statements are equivalent.
\begin{eqc}
\item $M$ is nilpotent, $\mathrm{fd}_{B}M<\infty$, $\mathrm{fd}M_{B}<\infty$, and $\mathrm{Tor}_{i}^{B}(M,M^{\oo j})=0$ for all $i,j\geq1$.
\item The functor $M\oo_{B}-: \BMod \ra \BMod$ is f-perfect and nilpotent.
\item The functor $-\oo_{B}M: \ModB \ra \ModB$ is f-perfect and nilpotent.
\end{eqc}
\end{prp}
\begin{prf*}
\proofofimp{i}{ii} By Lemma \ref{left perfect module}, $M\oo_{B}-$ is left perfect and nilpotent.
Since $\mathrm{fd}_{B}M<\infty$, it follows that $M\oo_{B}-$ maps projective modules to modules of finite flat dimension.
Hence the functor $M\oo_{B}-: \BMod \ra \BMod$ is f-perfect and nilpotent.

\proofofimp{ii}{i}
Assume that $M\otimes_{B}-:\BMod \ra \BMod$ is f-perfect and nilpotent.
Then it is easy to check that $M$ is nilpotent, and there exists $n\geq 1$ such that
$\mathrm{Tor}_{i}^{B}(M^{\oo j},-)=0$ for all $i,j\geq1$ with $i+j\geq n+1.$
Whenever $j=1$, we know that $\mathrm{Tor}_{i}^{B}(M,-)=0$ for $i\geq n,$
which yields that $\mathrm{fd}M_{B}<\infty.$
On the other hand, since $M\otimes_{B}-$ maps projective modules to the modules of finite flat dimension, it follows that $\mathrm{fd}_{B}M<\infty.$
Moreover, according to the condition (a) of Definition \ref{def:left perfect}, we have
$\mathrm{Tor}_{i}^{B}(M^{\oo j},M)=0$ for all $i,j\geq1$.
Note that $\mathrm{Tor}_{i}^{B}(M,M^{\oo j})=0$ for all $i,j\geq1$ if and only if $\mathrm{Tor}_{i}^{B}(M^{\oo j},M)=0$ for all $i,j\geq1$; see \cite[Lemma 5.2]{KP}. Then the implication $(ii)\Longrightarrow (i)$ holds.

The equivalence $(i)\Longleftrightarrow (iii)$ can be proved by
considering the opposite ring $B^\circ$ using the result \cite[Lemma 5.2]{KP}.
\end{prf*}



The next definition is from \cite{KP}.

\begin{dfn}
Let $\mathsf{F}: \BMod \ra \BMod$ be a left perfect right exact functor. A left $B$-module $X$ is called $\sf F$-projective if $\mathbb{L}_{i}\mathsf{F}^{j}(X)=0$ for all $i,j\geq1$.
\end{dfn}

\begin{rmk}
Let $\mathsf{F}: \BMod \ra \BMod$ be a left perfect right exact functor. It follows from \cite[Lemma 6.5]{KP} that a left $B$-module $X$ is $\sf F$-projective if and only if $\mathbb{L}_{i}\mathsf{F}(\mathsf{F}^{j}(X))=0$ for all $i\geq1$ and all $j\geq0$ if and only if $\mathbb{L}_{i}\mathsf{F}^{j}(\mathsf{F}^{s}(X))=0$ for all $i,j\geq1$ and all $s\geq0$.
\end{rmk}

In the following we collect some results on $\sf F$-projective modules, which can be found in \cite{Kostas}.
\begin{lem}\label{lem:properties of cleft extension}
Let $(\BMod,\AMod,\mathsf{i},\mathsf{e},\mathsf{l})$ be a cleft extension of module categories such that the induced endofunctor $\mathsf{F}$ is left perfect. Then the following statements hold:
\begin{prt}
\item For every projective left $A$-module $P$, $\mathsf{e}(P)$ is $\mathsf{F}$-projective.
\item If $X$ is an $\sf F$-projective left $B$-module, then so is $\mathsf{F}^{j}(X)$ for all $j\geq 1.$
\item If $X$ is an $\sf F$-projective left $B$-module, then there is an isomorphism $$\Ext_{A}^{i}(\mathsf{l}(X),Y)\cong \Ext_{B}^{i}(X,\mathsf{e}(Y))$$
    for all $i\geq1$ and every left $A$-module $Y$.
\end{prt}
\end{lem}

Let $R$ be a ring. A left $R$-module $C$ is called {\it cotorsion} if $\Ext_R^1(F,C)=0$ for each flat left $R$-module $F$. For a left $R$-module $M$, it is known that
$$\mathrm{fd}_{R}M=\sup\{i\in\mathbb{Z} \mid\Ext_{R}^{i}(M,C)\neq 0~\text{for some cotorsion left}~ R\text{-module}~C\}.$$

\begin{lem}\label{lem:e-flat dimension}
Let $(\BMod,\AMod,\mathsf{i},\mathsf{e},\mathsf{l})$ be a cleft extension of module categories such that the induced endofunctor $\mathsf{F}$ is  f-perfect.
Then $\mathrm{fd}_{B}\mathsf{e}(P)<\infty$ for all projective left $A$-module $P.$
\end{lem}
\begin{prf*}
Since every projective left $A$-module is a direct summand of $\mathsf{l}(Q)$ for some projective left $B$-module $Q$ by \cite[Lemma 2.3]{Kostas}, it enough to prove $\mathrm{fd}_{B}\mathsf{el}(Q)<\infty$.
Note that $\mathsf{el}(Q)\cong Q\oplus \mathsf{F}(Q)$ and $\mathrm{fd}_{B}(\mathsf{F}(Q))<\infty$ as $\sf F$ is f-perfect. It follows that $\mathrm{fd}_{B}\mathsf{el}(Q)<\infty$.
\end{prf*}

We end this section with the next result, which will be used later.

\begin{lem}\label{lem:refect flat dimension}
Let $(\BMod,\AMod,\mathsf{i},\mathsf{e},\mathsf{l})$ be a cleft extension of module categories such that the induced endofunctor $\mathsf{F}$ is f-perfect and nilpotent.
Assume that $\sf e$ preserves cotorsion modules and its right adjoint $\sf r$ is exact. Then for every $X\in \AMod$, one has $\mathrm{fd}_{A}X<\infty$ if and only if $\mathrm{fd}_{B}\mathsf{e}(X)<\infty$.
\end{lem}
\begin{prf*}
Keep the notation in the diagram (\ref{2.7.1}).

$(\Rightarrow)$: Let $X\in\AMod$ with $\mathrm{fd}_{A}X<\infty.$
Then there is an exact sequence $$0\ra Q_{k}\ra \cdots\ra Q_{1}\ra Q_{0}\ra X\ra 0$$
with $Q_{i}$ flat for $0\leq i\leq k$.
Applying $\sf e$ to the above exact sequence we get the following exact sequence
$$0\ra \mathsf{e}(Q_{k})\ra \cdots \ra\mathsf{e}(Q_{1})\ra \mathsf{e}(Q_{0})\ra \mathsf{e}(X)\ra 0.$$
Since $\sf r$ is exact and $\mathsf{e}$ preserves filtered colimits, it follows that $\mathsf{e}(Q_{i})$ is flat for $0\leq i\leq k$.
Hence $\mathrm{fd}_{B}\mathsf{e}(X)<\infty.$

$(\Leftarrow)$: Assume that $\mathrm{fd}_{B}\mathsf{e}(X)<\infty.$
Consider the following exact sequence
$$0\ra X'\ra P_{n-2}\ra \cdots \ra P_{1}\ra P_{0}\ra X\ra 0,$$
where $P_{i}$ is projective for $0\leq i\leq n-2$, and $n$ is a positive integer such that $\mathbb{L}_{p}\mathsf{F}^{q}=0$ for all $p,q\geq 1$ with $p+q\geq n+1$ as $\sf F$ is f-perfect; set $X'=X$ when $n=1$.
Applying exact functor $\sf e$ to the above exact sequence, we get the following exact sequence
$$0\ra \mathsf{e}(X')\ra \mathsf{e}(P_{n-2})\ra \cdots \ra \mathsf{e}(P_{1})\ra \mathsf{e}(P_{0})\ra \mathsf{e}(X)\ra 0.$$
By Lemma \ref{lem:properties of cleft extension}(a), we know that $\mathsf{e}(P_{n-2}),\cdots, \mathsf{e}(P_{0})$ are $\mathsf{F}$-projective.
Therefore by the dimension shifting, we have
$$\mathbb{L}_{i}\mathsf{F}^{j}(\mathsf e(X'))\cong \mathbb{L}_{i+n-1}\mathsf{F}^{j}(\mathsf e(X))=0$$
for all $i, j\geq1$, where the equality holds as $\mathbb{L}_{p}\mathsf{F}^{q}=0$ for all $p,q\geq 1$ with $p+q\geq n+1$.
Consequently, $\mathsf{e}(X')$ is $\mathsf{F}$-projective.
Since $\mathsf{F}$ is f-perfect, it follows from Lemma \ref{lem:e-flat dimension} that $\mathrm{fd}_{B}\mathsf{e}(P)<\infty$ for every projective left $A$-module $P$.
Thus $\mathrm{fd}_{B}\mathsf{e}(X')<\infty$.
Fix a flat resolution of $\mathsf{e}(X')$
$$0\ra Q_{t}\ra \cdots \ra Q_{1}\ra Q_{0}\ra \mathsf{e}(X')\ra 0$$
with $Q_t$ flat and $Q_i$ projective for $0\leq i\leq t-1$.
Since $\mathsf{e}(X')$ is $\mathsf{F}$-projective, we have $\mathbb{L}_{i}\mathsf{F}(\mathsf{e}(X'))=0$ for all $i\geq1$. Thus we get the following exact sequence
$$0\ra \mathsf{F}(Q_{t})\ra \cdots \ra \mathsf{F}(Q_{1})\ra \mathsf{F}(Q_{0})\ra \mathsf{F}\mathsf{e}(X')\ra 0.$$
Note that $\mathsf{er}\simeq \mathsf{F}'\oplus \mathrm{Id}$ and $\sf r$ is exact by assumption. It follows that $\mathsf{F}'$ is exact.
Then $\mathsf{F}$ preserves projective $B$-modules as $(\mathsf{F},\mathsf{F}')$ is an adjoint pair.
Since $\mathsf{F}$ preserves filtered colimits, $\mathsf{F}$ preserves flat $B$-modules. Thus ${\sf F}(Q_t)$ is flat and ${\sf F}(Q_i)$ is projective for $0\leq i\leq t-1$.
So $\mathrm{fd}_{B}\mathsf{F}\mathsf{e}(X')<\infty$.
Using an inductive argument, we get that $\mathrm{fd}_{B}\mathsf{F}^{j}\mathsf{e}(X')<\infty$ for all $j\geq 1$.
We mention that $\mathsf{F}^{j}\mathsf{e}(X')$ is $\sf F$-projective by Lemma \ref{lem:properties of cleft extension}(b).
It follows from Lemma \ref{lem:properties of cleft extension}(c) that $\Ext_{A}^{i}(\mathsf{l}\mathsf{F}^{j}\mathsf{e}(X'),V)\cong \Ext_{B}^{i}(\mathsf{F}^{j}\mathsf{e}(X'),\mathsf{e}(V))=0$
for each cotorsion left $A$-module $V$ and every $i\geq \mathrm{fd}_{B}\mathsf{F}^{j}\mathsf{e}(X')$, where the equality holds as $\mathsf{e}(V)$ is cotorsion by the assumption. This implies that $\mathrm{fd}_{A}\mathsf{l}\mathsf{F}^{j}\mathsf{e}(X')<\infty.$
We mention that $\mathsf{F}^{s}=0$ for some positive integer by assumption.
Then we consider the following exact sequences in $\AMod$:
\begin{equation*}
\begin{cases}
0\ra \mathsf{G}(X')\ra \mathsf{le}(X')\ra X'\ra 0;\\
0\ra  \mathsf{G}^{2}(X')\ra \mathsf{l}\mathsf{F}\mathsf{e}(X')\ra  \mathsf{G}(X')\ra 0;\\
\;\;\;\;\;\;\;\;\;\;\;\;\;\;\vdots\\
0\ra \mathsf{G}^{s-1}(X')\ra \mathsf{l}\mathsf{F}^{s-2}\mathsf{e}(X')\ra \mathsf{G}^{s-2}(X')\ra 0;\\
0\ra \mathsf{G}^{s}(X')\ra \mathsf{l}\mathsf{F}^{s-1}\mathsf{e}(X')\ra \mathsf{G}^{s-1}(X')\ra 0.
\end{cases}
\end{equation*}
It follows from Lemma \ref{lem of cleft extension} that $\mathsf{G}^{s}=0$  and so $\mathsf{G}^{s-1}(X')\cong \mathsf{l}\mathsf{F}^{s-1}\mathsf{e}(X')$.
Hence we get the exact sequence
$$0\ra\mathsf{l}\mathsf{F}^{s-1}\mathsf{e}(X') \ra\mathsf{l}\mathsf{F}^{s-2}\mathsf{e}(X')\ra \cdots\ra \mathsf{l}\mathsf{F}\mathsf{e}(X')\ra \mathsf{le}(X')\ra X'\ra 0.$$
We mention that $\mathsf{e}(X')$ is $\mathsf{F}$-projective and $\mathsf{e}$ preserves cotorsion modules. It follows from Lemma \ref{lem:properties of cleft extension}(c) that $\mathrm{fd}_{A}\mathsf{le}(X')<\infty.$
Thus one has $\mathrm{fd}_{A}X'<\infty$, and so $\mathrm{fd}_{A}X<\infty.$
\end{prf*}

\section{Gorenstein weak global dimension under cleft extensions}\label{gw}
\noindent
In this section, we study the homological behavior of the Gorenstein weak global dimension under cleft extensions of module categories.

\begin{ipg}
Let $R$ be a ring. An acyclic complex $T$ of flat left $R$-modules is called \emph{F-totally acyclic} if $E\otimes_R T$ is acyclic for every injective right $R$-module $E$.
A left $R$-module $G$ is called \emph{Gorenstein flat} if there exists an $F$-totally acyclic complex $T$ such that $G\cong\coker(T_1\to T_0)$. The Gorenstein flat dimension of a left $R$-module $M$, $\Gfd_RM$, is defined by declaring that $\Gfd_RM\leq n$ if and only if there is an exact sequence $0 \to G_n \to \cdots \to G_0 \to M \to 0$ of left $R$-modules with each $G_i$ Gorenstein flat. We denote by $\Gfd N_R$ the Gorenstein flat dimension of a right $R$-module $N$. We mention that the two values
$$\sup\{\Gfd_RM\ |\ M\ \text{is a left}\ R\text{-module}\}$$
and
$$\sup\{\Gfd N_R\ |\ N\ \text{is a right}\ R\text{-module}\}$$
are equal; see \cite[Corollary 1.4]{CET}. The \emph{Gorenstein weak global dimension} of $R$, denoted by $\Gw(R)$, is defined as the value
$$\sup\{\Gfd_RM\ |\ M\ \text{is a left}\ R\text{-module}\}=\sup\{\Gfd N_R\ |\ N\ \text{is a right}\ R\text{-module}\}.$$
\end{ipg}

\begin{ipg}
Let $R$ be a ring. We denote by $\mathrm{sfli}(R)$ (resp., $\mathrm{sfli}(\Rop)$) the supremum of the flat dimensions of injective left $R$-modules (resp., right $R$-modules), and denote by $\FFD(R)$ (resp., $\FFD(\Rop)$) the supremum of the flat dimensions of those left $R$-modules (resp., right $R$-modules) that have finite flat dimension.
\end{ipg}

By \cite[Theorem 2.8]{DLW} and \cite[Theorem 5.3]{Emmanouil}, we have the next result.

\begin{lem}\label{lem:Gw}
The following statements are equivalent for a nonnegative integer $n$.
\begin{eqc}
\item $\Gw(R)\leq n$.
\item $\mathrm{sfli}(R)\leq n$ and $\mathrm{sfli}(\Rop)\leq n$.
\item $\mathrm{sfli}(R)=\FFD(R)=\mathrm{sfli}(\Rop)=\FFD(\Rop)\leq n$.
\end{eqc}
If these equivalent conditions are satisfied, then there are equations
$$\Gw(R)=\mathrm{sfli}(R)=\FFD(R)=\mathrm{sfli}(\Rop)=\FFD(\Rop).$$
\end{lem}

\begin{lem}\label{lem:sfli}
Let $M$ be a left $R$-module. If $\mathrm{id}_{R}M<\infty$,
then $\mathrm{fd}_{R}M\leq \mathrm{sfli}(R).$
\end{lem}
\begin{prf*}
If $\mathrm{id}_{R}M=0,$ then the result holds by the definition of $\mathrm{sfli}(R).$
We proceed by induction on $\mathrm{id}_{R}M.$
If $\mathrm{id}_{R}M=n$, then there is an exact sequence
$$0\ra M\ra I\ra V\ra 0$$ with $I$ injective and
$\mathrm{id}_{R}V=n-1.$
Then $\fd_{R}I\leq \mathrm{sfli}(R)$ and $\fd_{R}V\leq \mathrm{sfli}(R)$ by induction hypothesis, and so $\mathrm{fd}_{R}M\leq \mathrm{sfli}(R).$
\end{prf*}

\begin{lem}\label{lem:preserves cotorsion module}
Let ${\sf l}: \BMod\to\AMod$ and ${\sf e}: \AMod\to\BMod$ be two exact functors such that $(\sf l,\sf e)$ is an adjoint pair. Then $\sf e$ preserves cotorsion modules.
\end{lem}
\begin{prf*}
Fix a cotorsion left $A$-module $C$ and a flat left $B$-module $X.$
Since $\sf l$ and $\sf e$ are exact, one gets an isomorphism $$\Ext^{1}_{B}(X,\mathsf{e}(C))\cong \Ext^{1}_{A}(\mathsf{l}(X),C).$$
Since $\sf l$ preserves projective modules and filtered colimits, $\sf l$ preserves flat modules.
It follows that $\Ext^{1}_{A}(\mathsf{l}(X),C)=0$, which implies that $\mathsf{e}(C)$ is a cotorsion left $B$-module.
\end{prf*}

\begin{lem}\label{prop:sfli}
Let $(\BMod,\AMod,\mathsf{i},\mathsf{e},\mathsf{l})$ be a cleft extension of module categories. Assume that the induced endofunctor $\mathsf{F}$ on $\BMod$ is nilpotent with $s$ an integer such that $\mathsf{F}^{s}=0$. Then the following statements hold:
\begin{prt}
\item If $\sf r$ is exact, then there is an inequality $\mathrm{sfli}(B)\leq \mathrm{sfli}(A).$
\item If $\mathsf{F}$ is f-perfect and both $\mathsf{r}$ and $\mathsf{l}$ are exact,
then there is an inequality
$$\mathrm{sfli}(A)\leq \mathrm{sfli}(B)+s-1.$$
\end{prt}
\end{lem}
\begin{prf*}
Keep the notation in the diagram (\ref{2.7.1}).

(a) We may assume that $\mathrm{sfli}(A)<\infty$. Let $V$ be a cotorsion left $B$-module. Since $\sf r$ is exact, we get the following isomorphism for any injective left $B$-module $I$ and every $k\geq1$:
\begin{align*}
\Ext_{A}^{k}(\mathsf{r}(I),\mathsf{r}(V))&\cong\Ext_{B}^{k}(\mathsf{er}(I),V)\\
&\cong \Ext_{B}^{k}(I,V)\oplus \Ext_{B}^{k}(\mathsf{F}'(I),V),
\end{align*}
where the second isomorphism holds as $\mathsf{er}\simeq \mathrm{Id}\oplus \mathsf{F}'.$
From Lemma \ref{lem:preserves cotorsion module} we know that $\mathsf{r}(V)$ is a cotorsion left $A$-module.
Then we have $\Ext_{A}^{k}(\mathsf{r}(I),\mathsf{r}(V))=0$ for all $k\geq \mathrm{sfli}(A)+1$ as $\mathsf{r}(I)$ is injective.
It implies that $\Ext_{B}^{k}(I,V)=0$ for all $k\geq \mathrm{sfli}(A)+1,$
which yields that  $\mathrm{sfli}(B)\leq \mathrm{sfli}(A).$

(b) We assume that $\mathrm{sfli}(B)<\infty$. By \cite[Lemma 2.6]{Kostas} that every injective module in $\AMod$ is a direct summands of $\mathsf{r}(I)$ for some injective left $B$-module $I.$
Let $X$ be a cotorsion left $A$-module. Since $\sf r$ is exact, we have the following isomorphism
\begin{align*}
\Ext_{A}^{k}(\mathsf{r}(I),\mathsf{r}\mathsf{F}'^{j}\mathsf{e}(X))&\cong\Ext_{B}^{k}(\mathsf{er}(I),\mathsf{F}'^{j}\mathsf{e}(X))\\
&\cong \Ext_{B}^{k}(I,\mathsf{F}'^{j}\mathsf{e}(X))\oplus \Ext_{B}^{k}(\mathsf{F}'(I),\mathsf{F}'^{j}\mathsf{e}(X))
\end{align*}
for $j\geq 0$,
where the second isomorphism holds as $\mathsf{er}\simeq \mathrm{Id} \oplus \mathsf{F}'.$
 Since $\sf F$ is f-perfect, one has $\mathsf{F}\simeq M\oo_{B}-$ and $\mathsf{F}'\simeq \Hom_{B}(M,-)$ for some $B$-bimodule $M$ with $\mathrm{fd}M_{B}<\infty$; see Remark \ref{remark:fg} and Proposition \ref{perfect module}.
 Thus $\id_{B}\mathsf{F}'(I)=\id_{B}\Hom_{B}(M,I)<\infty$ for any injective $B$-module $I$, and we know from Lemma \ref{lem:sfli} that $\fd_{B}I\leq\mathrm{sfli}(B)$ and $\fd_{B}\mathsf{F}'(I)\leq\mathrm{sfli}(B).$
We mention that $\mathsf{el}\simeq \mathsf{F}\oplus \mathrm{Id}$ and $\mathsf{er}\simeq \mathsf{F}'\oplus \mathrm{Id}$.
Since $\sf l$ and $\sf r$ are exact by assumption, so are $\mathsf{F}$ and $\mathsf{F}'$.
Note that $(\mathsf{F},\mathsf{F}')$ is an adjoint pair.
It follows from Lemma \ref{lem:preserves cotorsion module} that
 $\mathsf{F}'^{j}\mathsf{e}(X)$ is cotorsion left $B$-module for $j\geq0$.
Thus $\Ext_{B}^{k}(I,\mathsf{F}'^{j}\mathsf{e}(X))=0= \Ext_{B}^{k}(\mathsf{F}'(I),\mathsf{F}'^{j}\mathsf{e}(X))$ for all $k\geq \mathrm{sfli}(B)+1$, and so $\Ext_{A}^{k}(\mathsf{r}(I),\mathsf{r}\mathsf{F}'^{j}\mathsf{e}(X))=0$
for all $k\geq \mathrm{sfli}(B)+1.$
Consider the following exact sequences:
\begin{equation*}
\begin{cases}
0\ra X\ra \mathsf{re}(X)\ra \mathsf{G}'(X)\ra 0;\\
0\ra \mathsf{G}'(X)\ra \mathsf{r}\mathsf{F}'\mathsf{e}(X)\ra \mathsf{G}'^{2}(X)\ra 0;\\
\;\;\;\;\;\;\;\;\;\;\;\;\;\;\vdots\\
0\ra \mathsf{G}'^{s-2}(X)\ra \mathsf{r}\mathsf{F}'^{s-2}\mathsf{e}(X)\ra \mathsf{G}'^{s-1}(X)\ra 0.\\
0\ra \mathsf{G}'^{s-1}(X)\ra \mathsf{r}\mathsf{F}'^{s-1}\mathsf{e}(X)\ra \mathsf{G}'^{s}(X)\ra 0.\\
\end{cases}
\end{equation*}
Note that $\mathsf{F}^{s}=0$.
It follows from \cite[Proposition 4.19]{Kostas} that $\mathsf{F}'^{s}=0.$
Thus we have $\mathsf{G}'^{s}=0$ by Lemma \ref{lem of cleft coextension}(a) and so $\mathsf{G}'^{s-1}(X)\cong \mathsf{r}\mathsf{F}'^{s-1}\mathsf{e}(X)$.
Therefore there is an exact sequence of $\AMod:$
$$0\ra X\ra \mathsf{re}(X) \ra \mathsf{r}\mathsf{F}'\mathsf{e}(X)\ra \cdots\ra \mathsf{r}\mathsf{F}'^{s-1}\mathsf{e}(X)\ra  0.$$
It follows that
\begin{align*}
\Ext_{A}^{k+s-1}(\mathsf{r}(I),X)
\cong \Ext_{A}^{k}(\mathsf{r}(I),\mathsf{r}\mathsf{F}'^{s-1}\mathsf{e}(X))
=0
\end{align*}
for all $k\geq \mathrm{sfli}(B)+1.$
Consequently, $\mathrm{sfli}(A)\leq \mathrm{sfli}(B)+s-1.$
\end{prf*}

We mention that if $({\BMod},{\AMod},{\sf i}_{1},\mathsf{e}_{1},\mathsf{l}_{1})$ is a cleft extension of module categories, then there exists a cleft extension $(\ModB,\ModA,\mathsf{i}_{2},\mathsf{e}_{2},\mathsf{l}_{2})$; see Remark \ref{remark:fg}.
In the following, we call $(\ModB,\ModA,\mathsf{i}_{2},\mathsf{e}_{2},\mathsf{l}_{2})$ the induced cleft extension by $({\BMod},{\AMod},{\sf i}_{1},\mathsf{e}_{1},\mathsf{l}_{1})$.
Let $\mathsf{F}_{1}$ and $\mathsf{F}_{2}$ be the induced endofunctors on $\BMod$ and $\ModB$ respectively, and denote by $\mathsf{r}_{1}$ the right adjoint of $\mathsf{e}_{1}$, and by $\mathsf{r}_{2}$ the right adjoint of $\mathsf{e}_{2}$. We then give the main result in this section.

\begin{thm}\label{thm:Gw}
Let $(\BMod,\AMod,\mathsf{i}_{1},\mathsf{e}_{1},\mathsf{l}_{1})$ be a cleft extension of module categories and $(\ModB,\ModA,\mathsf{i}_{2},\mathsf{e}_{2},\mathsf{l}_{2})$ be the induced cleft extension. Assume that the induced endofunctor $\mathsf{F}_{1}$ on $\BMod$ is f-perfect and nilpotent with $s$ an integer such that $\mathsf{F}_{1}^{s}=0$. If $\mathsf{r}_{1}$ and $\mathsf{r}_{2}$ are exact,
then there are inequalities
$$\Gw (B)\leq \Gw (A) \leq \Gw (B)+s-1,$$
In particular, $\Gw(B)$ is finite if and only if $\Gw(A)$ is finite.
\end{thm}
\begin{prf*}
Since $\mathsf{r}_{1}$ and $\mathsf{r}_{2}$ are exact, so are $\mathsf{l}_{1}$ and $\mathsf{l}_{2}$ by Lemma \ref{lem:left-right}.
According to the assumption that $\mathsf{F}_{1}$ is f-perfect and $\mathsf{F}_{1}^{s}=0$, it follows from Remark \ref{remark:fg} and Proposition \ref{perfect module} that $\mathsf{F}_{2}$ is f-perfect and $\mathsf{F}_{2}^{s}=0$.
Then the inequalities in the statement follow from Lemmas \ref{lem:Gw} and \ref{prop:sfli}.
\end{prf*}

\begin{rmk}
To push the Gorenstein flat dimension beyond the setting of coherent rings, Christensen, Estrada, Liang, Thompson, Wu and Yang \cite{CELTWY} introduced the Gorenstein flat-cotorsion dimension, which refines the Gorenstein flat dimension and exhibits properties more aligned with expectations. The corresponding global dimension is termed the Gorenstein flat-cotorsion global dimension, whose finiteness is explored in \cite{CELTW}. A recent result by Weiqing Li shows that the Gorenstein weak global dimension $\Gw(R)$ and Gorenstein flat-cotorsion global dimension $\Gfc(R)$ are the same for arbitrary associative ring $R$; see \cite[Theorem 2.3]{Li}. Thus under the same assumptions as in Theorem \ref{thm:Gw}, one gets the inequalities
$\Gfc(B)\leq \Gfc(A) \leq \Gfc(B)+s-1$.
\end{rmk}

\section{Applications to some known ring extensions}

\noindent
In this section, we apply the main results in Section \ref{gw} on the Gorenstein weak global dimension to $\theta$-extensions and Morita context rings; one refers to Examples \ref{ex:extension} and \ref{ex:Morita ring} for the definitions of $\theta$-extensions and Morita context rings, respectively.

\begin{thm}\label{IT:corollary}
Let $M$ be a $B$-bimodule with $M^{\oo s}=0$ for some positive integer $s$, and let $\theta:M\oo_{B}M\ra M$ be an associative $B$-bimodule homomorphism. If $M$ is projective as a left $B$-module and a right $B$-module, then there are inequalities
    $$\Gw(B)\leq\Gw(\G)\leq \Gw(B)+s-1.$$
    In particular, $\Gw(B)$ is finite if and only if $\Gw(\G)$ is finite.
\end{thm}
\begin{prf*}
By Example \ref{ex:extension}, we have the following cleft extensions
\begin{equation*}
  \xymatrix@C=2pc{
    \BMod \ar[r]^-{\mathsf{Z}_{1}} & \eMod \ar[r]^-{\mathsf{U}_{1}} & \BMod
    \ar@/_1.8pc/[l]_-{\mathsf{T}_{1}=(B\ltimes_{\theta}M)\oo_{B}-}
    \ar@/^1.8pc/[l]^-{\mathsf{H}_{1}=\mathrm{Hom}_{B}(B\ltimes_{\theta}M,-)}\ar@(lu,ru)[]^{\sf F_1}}
\end{equation*}
and
\begin{equation*}
  \xymatrix@C=2pc{
    \ModB \ar[r]^-{\mathsf{Z}_{2}} & \Mode \ar[r]^-{\mathsf{U}_{2}} & \ModB.
    \ar@/_1.8pc/[l]_-{\mathsf{T}_{2}=-\oo_{B}(B\ltimes_{\theta}M)}
    \ar@/^1.8pc/[l]^-{\mathsf{H}_{2}=\mathrm{Hom}_{B}(B\ltimes_{\theta}M,-)}\ar@(lu,ru)[]^{\sf F_2}}
\end{equation*}
with $\mathsf{F}_{1}\simeq M\oo_{B}-$ and $\mathsf{F}_{2}\simeq -\oo_{B}M.$

Since $_{B}M$ and $M_{B}$ are projective, one gets that $\mathsf{F}_{1}$ and $\mathsf{F}_{2}$ are f-perfect, and $\mathsf{H}_{1}$ and $\mathsf{H}_{2}$ are exact. Hence the statement follows from Theorem \ref{thm:Gw}.
\end{prf*}

Let $M$ be a $B$-bimodule. We mention that the tensor ring $T_{B}(M)$ can be viewed as $\theta$-extension: if we set $M'=M\oplus M^{\oo 2}\oplus\cdots$, then there is an isomorphism $T_{B}(M)\cong B\ltimes_{\theta}M'$ with the $B$-bimodule homomorphism $\theta: M'\otimes_B M'\to M'$ induced by the natural homomorphisms $M^{\oo k}\otimes_B M^{\oo l} \to M^{\oo (k+l)}$; see \cite[Example 6.8]{KP}. Thus the next result is an immediate consequence of Theorem \ref{IT:corollary}.

\begin{cor}\label{new}
Let $M$ be a $B$-bimodule with $M^{\oo s}=0$ for some positive integer $s$. If $M$ is projective as a left $B$-module and a right $B$-module, then there are inequalities
    $$\Gw(B)\leq\Gw(T_{B}(M))\leq \Gw(B)+s-1.$$
    In particular, $\Gw(B)$ is finite if and only if $\Gw(T_{B}(M))$ is finite.
\end{cor}

We then apply Theorem \ref{IT:corollary} to Morita context rings.

\begin{cor}\label{thm:Morita ring}
Let $\Lambda =\left(\begin{smallmatrix}  A & {_{A}}N_{B}\\  {_{B}}M_{A} & B \\\end{smallmatrix}\right)$
be a Morita context ring with $N\oo_{B}M=0=M\oo_{A}N$.
If $M$ is projective as a right $A$-module and a left $B$-module, and $N$ is projective as a left $A$-module and a right $B$-module, then there are inequalities
$$\max\{\Gw (A),\Gw (B)\}\leq\Gw (\Lambda)$$
and
$$\Gw (\Lambda)\leq \max\{\Gw (A),\Gw (B)\}+1.$$
In particular, $\Gw(\Lambda)$ is finite if and only if the invariants $\Gw(A)$ and $\Gw(B)$ are finite.
\end{cor}
\begin{prf*}
We mention that $M\oplus N$ is an $A\times B$-bimodule and $\Lambda\simeq(A\times B)\ltimes (M\oplus N)$; see \cite[Proposition 2.5]{GP}. Since $M_{A}, N_{B}, {_{A}N}$ and $_{B}M$ are projective, it is easy to check that $M\oplus N$ is projective as a left $A\times B$-module and a right $A\times B$-module. On the other hand, one gets that $(M\oplus N)^{\oo 2}=0$ as $N\oo_{B}M=0=M\oo_{A}N$. Thus the statement follows from Theorem \ref{IT:corollary}.
\end{prf*}

\section{Singularity categories with respect to flat-cotorsion modules}
\noindent
Let $\mathsf{A}$ be an abelian subcategory and $\mathsf{W}$ a full subcategory of $\sf A$.
If $\mathsf{W}$ is self-orthogonal, that is $\Ext_{\mathsf{A}}^{n}(X,Y)=0$ for any $X,Y\in \mathsf{W}$ and $n\geq 1,$ then Chen \cite{Chen} defined a relative singularity category $\mathbf{D}_{\mathsf{W}}(\mathsf{A}):=\mathbf{D}^{b}(\mathsf{A})/\mathbf{K}^{b}(\mathsf{W})$.
In what follows, for a ring $R$ we write $\mathbf{D}^{b}(R)$ (resp., $\mathbf{D}(R)$) for $\mathbf{D}^{b}(\RMod)$ (resp., $\mathbf{D}(\RMod)$) and denote by $\mathsf{FlatCot}(R)$ the subcategory consisting of flat-cotorsion $R$-modules; it is self-orthogonal.
In this section, we consider the relative singularity category
$\mathbf{D}_{\mathsf{FlatCot}}(R): =\mathbf{D}^{b}(R)/\mathbf{K}^{b}(\mathsf{FlatCot}(R))$ under the context of cleft extensions of module categories.
In the following, we will investigate relative singularity categories with respect to flat-cotorsion modules in terms of flat model structures on category of chain complexes constructed by Gillespie \cite{Gillespie}.

\begin{ipg}
We let $\mathsf{Ch}(R)$ denote the category of complexes of left $R$-modules (or $R$-complexes). Recall from \cite{Gillespie} that an $R$-complex $X$ is called {\it dg-flat} if  each component $X^{n}$ is flat and $\Hom_R(X,C)$ is exact whenever $C$ is an acyclic $R$-complex with each cycles $\mathrm{Z}^{n}C$ cotorsion. An $R$-complex $X$ is called {\it dg-cotorsion} if each component $X^{n}$ is cotorsion and $\Hom_R(F,X)$ is exact whenever $F$ is an acyclic $R$-complex with each cycles $\mathrm{Z}^{n}F$ flat. We denote by $\dg\widetilde{\mathsf{F}}$ (resp., $\dg\widetilde{\mathsf{C}}$, $\mathsf{Acy}$) the class of dg-flat (resp., dg-cotorsion, acyclic) $R$-complexes. We mention that every complex of cotorsion $R$-modules is dg-cotorsion; see \cite[Theorem 1.3]{Ba}.
\end{ipg}

The next result can be found in \cite[Corollary 5.1]{Gillespie}.

\begin{lem}\label{lem:flat model structure}
There is a flat model structure $(\dg\widetilde{\mathsf{F}},\mathsf{Acy},\dg\widetilde{\mathsf{C}})$ on $\mathsf{Ch}(R)$; its homotopy category is $\mathbf{D}(R).$
\end{lem}

The proof of the following result is routine, which can be found in \cite[Lemma 4.5]{LM}.

\begin{lem}\label{quasi-iso}
Let $X$ be a bounded $R$-complex with  each component $X^{n}$ cotorsion and $\mathrm{fd}_{R}X^{n}<\infty$.
Then $X$ is quasi-isomorphic to a bounded complex in $\mathsf{FlatCot}(R)$.
\end{lem}

\begin{prp} \label{prop:relative singularity categories}
Let $A$ and $B$ be right coherent rings and $(\BMod,\AMod,\mathsf{i},\mathsf{e},\mathsf{l})$ a cleft extension of module categories. Assume that the following conditions are satisfied:
\begin{prt}
\item $\sf l$ is exact and preserves products;
\item $\mathbb{ L}_{i}\mathsf{q}=0$ for $i$ large enough, and $\sf q$ preserves products;
\item $\sf i$ preserves cotorsion modules and $\mathrm{fd}_{A}\mathsf{i}(V)<\infty$ for any flat-cotorsion $B$-module $V$;
\item $\mathrm{fd}_{B}\mathsf{e}(W)<\infty$ for any flat-cotorsion $A$-module $W$.
\end{prt}
Then there is a diagram
\begin{equation*}
  \xymatrix@C=2pc{
   \mathbf{D}_{\mathsf{FlatCot}}(B)\ar[r]^{\overline{\sf Ri}} & \mathbf{D}_{\mathsf{FlatCot}}(A) \ar@/_1.8pc/[l]_{\overline{\mathsf{Lq}}}
   \ar[r]^{\overline{\sf Re}} & \mathbf{D}_{\mathsf{FlatCot}}(B) \ar@/_1.8pc/[l]_{\overline{\mathsf{Ll}}}}
\end{equation*}
of relative singularity categories and triangle functors such that $(\overline{\sf Lq},\overline{\sf Ri})$ and $(\sf \overline{Ll},\overline{Re})$ are adjoint pairs, $\overline{\sf Re}\circ\overline{\sf Ri}\simeq \mathrm{Id}_{\mathbf{D}_{\mathsf{FlatCot}}(B)}$ and $\overline{\sf Lq}\circ\overline{\sf Ll}\simeq \mathrm{Id}_{\mathbf{D}_{\mathsf{FlatCot}}(B)}$.
\end{prp}
\begin{prf*}
We mention that a cleft extension of module categories
\begin{equation*}
  \xymatrix@C=2pc{
   \BMod \ar[r]^-{\sf i} & \AMod \ar@/_1.8pc/[l]_-{\sf q}
   \ar[r]^-{\sf e} & \BMod \ar@/_1.8pc/[l]_-{\sf l}}
\end{equation*}
induces a cleft extension of complex categories
\begin{equation*}
  \xymatrix@C=2pc{
   \mathsf{Ch}(B) \ar[r]^-{\sf i} & \mathsf{Ch}(A) \ar@/_1.8pc/[l]_{\sf q}
   \ar[r]^{\sf e} & \mathsf{Ch}(B). \ar@/_1.8pc/[l]_{\sf l}}
\end{equation*}
By Lemma \ref{lem:flat model structure}, there exist flat model structures on $\mathsf{Ch}(A)$ and $\mathsf{Ch}(B).$  We claim that $(\sf l,\sf e)$ and $(\sf q,\sf i)$ are Quillen adjunctions.
Since $\sf l$ is exact and $\sf i$ preserves cotorsion modules, it follows from Lemma \ref{lem:preserves cotorsion module} that $\sf i$ and $\sf e$ preserve dg-cotorsion complexes. We mention that $\sf i$ and $\sf e$ are exact. Then $\sf i$ and $\sf e$ preserve fibrations and trivial fibrations. It follows from \cite[Lemma 1.3.4]{Hovey}  that $(\sf l,\sf e)$ and $(\sf q,\sf i)$ are Quillen adjunctions.
By \cite[Lemma 1.3.10]{Hovey} we get the following diagram of derived categories
\begin{equation*}
  \xymatrix@C=2pc{
   \mathbf{D}(B) \ar[r]^-{\sf Ri} & \mathbf{D}(A) \ar@/_1.8pc/[l]_{\sf Lq}
   \ar[r]^{\sf Re} & \mathbf{D}(B) \ar@/_1.8pc/[l]_{\sf Ll}}
\end{equation*}
such that $(\sf Ll,Re)$ and $(\sf Lq,Ri)$ are adjoint pairs of triangle functors.
We claim that the above functors restricts to
\begin{equation*}
  \xymatrix@C=2pc{
   \mathbf{D}^{b}(B) \ar[r]^-{\sf Ri} & \mathbf{D}^{b}(A) \ar@/_1.8pc/[l]_{\sf Lq}
   \ar[r]^{\sf Re} & \mathbf{D}^{b}(B). \ar@/_1.8pc/[l]_{\sf Ll}}
\end{equation*}
For $X\in \mathbf{D}^{b}(A)$, by Lemma \ref{lem:flat model structure}, there exists an exact sequence of $A$-complexes
$$0\ra C\ra Q\ra X\ra 0$$
such that $Q$ is dg-flat and $C$ is acyclic dg-cotorsion.
In particular, $Q$ is quasi-isomorphic to $X$.
We may write
$$X=\cdots \ra X^{-2}\xrightarrow{d_{X}^{-2}} X^{-1}\xrightarrow{d_{X}^{-1}} X^{0}\ra 0\ra 0\ra \cdots$$
and $$Q=\cdots \ra Q^{-2}\xrightarrow{d_{Q}^{-2}} Q^{-1}\xrightarrow{d_{Q}^{-1}} Q^{0}\xrightarrow{d_{Q}^{0}} Q^{1}\xrightarrow{d_{Q}^{1}} Q^{2}\ra \cdots$$
Note that $Q^{i}$ is flat and $\mathsf{q}\simeq \mathsf{q}(A)\otimes_{A}-$ by the Eilenberg-Watts theorem. Then one has
$$\mathsf{H}^{i}(\mathsf{Lq}(X))=\mathsf{H}^{i}(\mathsf{q}(Q))\cong \mathbb{L}_{j-i}\mathsf{q}(\mathsf{B}^{j+1}(Q))$$
for all $i>0$ and all $j>i.$
According to the assumption that $\mathbb{L}_{i}\mathsf{q}=0$ for $i$ large enough, we have $\mathsf{H}^{i}(\mathsf{Lq}(X))=0$ for all $i>0.$
Since $X$ is a bounded complex, we may assume that $X^{-t}=0$ for any $t\geq m+1$ for some integer $m.$
Thus we have an exact sequence
$$\cdots \ra Q^{-(m+2)}\ra Q^{-(m+1)}\ra \ker d_{Q}^{-m}\ra 0.$$
It follows that $$\mathsf{H}^{-(m+i+1)}(\mathsf{Lq}(X))=\mathsf{H}^{-(m+i+1)}(\mathsf{q}(Q))\cong \mathbb{L}_{i}\mathsf{q}(\ker d_{Q}^{-m})=0$$ for enough large $i.$
So $\mathsf{Lq}$ can be restricted to a functor $\mathbf{D}^{b}(A)\ra \mathbf{D}^{b}(B)$.
We mention that $\sf l, e$ and $\sf i$ are exact.
Then it is easy to see that $\mathsf{Ll}, \mathsf{Re}$ and $\mathsf{Ri}$ can be restricted to bounded derived categories.
By conditions (a) and (b), it follows from Lemma \ref{lem:preserves cotorsion module} that $$\mathsf{Lq}(\mathbf{K}^{b}(\mathsf{FlatCot}(A)))\subseteq \mathbf{K}^{b}(\mathsf{FlatCot}(B))$$ and
$$\mathsf{Ll}(\mathbf{K}^{b}(\mathsf{FlatCot}(B)))\subseteq \mathbf{K}^{b}(\mathsf{FlatCot}(A)).$$
Since $\sf l$ is exact, it follows from Lemma \ref{lem:preserves cotorsion module} that $\sf e$ preserves cotorsion modules.
 By combining Lemma \ref{quasi-iso} and conditions (c) and (d) one deduces that $$\mathsf{Ri}(\mathbf{K}^{b}(\mathsf{FlatCot}(B)))\subseteq \mathbf{K}^{b}(\mathsf{FlatCot}(A))$$ and
$$\mathsf{Re}(\mathbf{K}^{b}(\mathsf{FlatCot}(A)))\subseteq \mathbf{K}^{b}(\mathsf{FlatCot}(B)).$$
Thus by \cite[Lemma 1.2]{Orlov} we get the following diagram of relative singularity categories
\begin{equation*}
  \xymatrix@C=2pc{
   \mathbf{D}_{\mathsf{FlatCot}}(B)\ar[r]^{\overline{\sf Ri}} & \mathbf{D}_{\mathsf{FlatCot}}(A) \ar@/_1.8pc/[l]_{\overline{\mathsf{Lq}}}
   \ar[r]^{\overline{\sf Re}} & \mathbf{D}_{\mathsf{FlatCot}}(B) \ar@/_1.8pc/[l]_{\overline{\mathsf{Ll}}}}
\end{equation*}
such that $(\overline{\sf Lq},\overline{\sf Ri})$ and $(\sf \overline{Ll},\overline{Re})$ are adjoint pairs, $\overline{\sf Re}\circ\overline{\sf Ri}\simeq \mathrm{Id}_{\mathbf{D}_{\mathsf{FlatCot}}(B)}$ and $\overline{\sf Lq}\circ\overline{\sf Ll}\simeq \mathrm{Id}_{\mathbf{D}_{\mathsf{FlatCot}}(B)}$.
\end{prf*}

\begin{thm}\label{thm:relative singularity categories}
Let $A$ and $B$ be right coherent rings and $(\BMod,\AMod,\mathsf{i},\mathsf{e},\mathsf{l})$ a cleft extension of module categories such that the induced endofunctor $\mathsf{F}$ is f-perfect and nilpotent. Assume that the following conditions are satisfied:
\begin{prt}
\item $\sf l$  and $\sf r$ are exact;
\item $\sf l$ and $\sf q$ preserve products;
\item $\sf i$ preserves cotorsion modules.
\end{prt}
Then there is a diagram
\begin{equation*}
  \xymatrix@C=2pc{
   \mathbf{D}_{\mathsf{FlatCot}}(B)\ar[r]^{\overline{\sf Ri}} & \mathbf{D}_{\mathsf{FlatCot}}(A) \ar@/_1.8pc/[l]_{\overline{\mathsf{Lq}}}
   \ar[r]^{\overline{\sf Re}} & \mathbf{D}_{\mathsf{FlatCot}}(B) \ar@/_1.8pc/[l]_{\overline{\mathsf{Ll}}}}
\end{equation*}
of  relative singularity categories and triangle functors such that $(\overline{\sf Lq},\overline{\sf Ri})$ and $(\overline{\sf Ll},\overline{\sf Re})$ are adjoint pairs, $\overline{\sf Re}\circ\overline{\sf Ri}\simeq \mathrm{Id}_{\mathbf{D}_{\mathsf{FlatCot}}(B)}$ and $\overline{\sf Lq}\circ\overline{\sf Ll}\simeq \mathrm{Id}_{\mathbf{D}_{\mathsf{FlatCot}}(B)}$.
\end{thm}
\begin{prf*}
Since $\sf F$ is f-perfect and nilpotent, it follows from \cite[Proposition 6.6]{KP} that $\mathbb{L}_{i}\mathsf{q}=0$ for $i$ large enough. Note that $\sf l$ is exact. We know from Lemma \ref{lem:preserves cotorsion module} that $\sf e$ preserves cotorsion modules.
 Consider a flat-cotorsion left $B$-module $V$.
Since $\sf r$ is exact, it follows from
Lemma \ref{lem:refect flat dimension} that $\mathrm{fd}_{A}\mathsf{i}(V)<\infty$ if and only if $\mathrm{fd}_{B}\mathsf{ei}(V)<\infty$. Since $\mathsf{ei}(V)\cong V$, we get $\mathrm{fd}_{A}\mathsf{i}(V)<\infty$. Therefore the result follows by Proposition \ref{prop:relative singularity categories}.
\end{prf*}

\begin{cor}\label{cor6.6}
Let $\Lambda =\left(\begin{smallmatrix}  A & {_{A}}N_{B}\\  {_{B}}M_{A} & B \\\end{smallmatrix}\right)$
be a Morita context ring with $A$ and $B$ right coherent and $N\oo_{B}M=0=M\oo_{A}N$.
Assume that $M$ is projective as a right $A$-module and a left $B$-module, and $N$ is projective as a left $A$-module and a right $B$-module, and assume that $M$ and $N$ are finitely presented as a right $A$-module and a right $B$-module respectively. Then there is a diagram
\begin{equation*}
  \xymatrix@C=2pc{
   \mathbf{D}_{\mathsf{FlatCot}}(A)\times\mathbf{D}_{\mathsf{FlatCot}}(B) \ar[r]^-{\overline{\sf Ri}} & \mathbf{D}_{\mathsf{FlatCot}}(\Lambda) \ar@/_1.8pc/[l]_{\overline{\mathsf{Lq}}}
   \ar[r]^-{\overline{\sf Re}} & \mathbf{D}_{\mathsf{FlatCot}}(A)\times\mathbf{D}_{\mathsf{FlatCot}}(B) \ar@/_1.8pc/[l]_-{\overline{\mathsf{Ll}}}}
\end{equation*}
of singularity categories and triangle functors such that $(\overline{\sf Lq},\overline{\sf Ri})$ and $(\sf \overline{Ll},\overline{Re})$ are adjoint pairs, $\overline{\sf Re}\circ\overline{\sf Ri}\simeq \mathrm{Id}_{\mathbf{D}_{\mathsf{FlatCot}}(A)\times\mathbf{D}_{\mathsf{FlatCot}}(B)}$ and $ \overline{\sf Lq}\circ\overline{\sf Ll}\simeq \mathrm{Id}_{\mathbf{D}_{\mathsf{FlatCot}}(A)\times\mathbf{D}_{\mathsf{FlatCot}}(B)}$.
\end{cor}
\begin{prf*}
By Example \ref{ex:Morita ring} we have a cleft extension of module categories
\begin{equation*}
  \xymatrix@C=2pc{
   \ABMod \ar[r]^-{\sf i} & \LMod \ar@/_1.8pc/[l]_-{\sf q}
   \ar[r]^-{\sf e} & \ABMod \ar@/_1.8pc/[l]_-{\sf l}}
\end{equation*}
Here $\mathsf{q}(X,Y,f,g)=(X,\coker f)\oplus (\coker g, Y)$,
$\mathsf{l}(X,Y)=(X,M\oo_{A}X,1,0)\oplus (N\oo_{B}Y,Y,0,1)$ and
$\mathsf{i}(X,Y)=(X,Y,0,0).$
Since $M_{A}$ and $N_{B}$ are finitely presented, it is easy to check that $\sf l$ and $\sf q$ preserve products.
By \cite[Theorem 4.4]{CRZ} and \cite[Remark 2.6]{LMY}, we know that $\sf i$ preserves cotorsion modules.
Since $M_{A}, N_{B}, {_{A}N}$ and $_{B}N$ are projective, it follows from Example
\ref{ex:Morita ring} that $\sf l$ and $\sf r$ are exact.
Consequently, we get the conclusion by Theorem \ref{thm:relative singularity categories}.
\end{prf*}

\section*{Acknowledgment}

\noindent
We thank Sergio Estrada for pointing out the fact that the Gorenstein weak
global dimension and the Gorenstein flat-cotorsion global dimension are the same for arbitrary associative ring, which was proved by Weiqing Li; see \cite[Theorem 2.3]{Li}.

L. Liang was partly supported by NSF of China grant 12271230. Y. Ma was partly supported by NSF of Gansu Province grant 23JRRA866, and the Youth Foundation of Lanzhou Jiaotong University grant 2023023. G. Yang was partly supported by NSF of China grant 12161049. All three authors were also partly supported by the Foundation for Innovative Fundamental Research Group Project of Gansu Province grants 23JRRA684 and 25JRRA805.


\bibliographystyle{amsplain-nodash}

\end{document}